\documentclass[12pt]{article}
\usepackage{anysize}
\usepackage{amssymb}
\usepackage{amsmath}
\usepackage{amsthm}
\usepackage{graphics}
\usepackage{epsfig}
\usepackage{url}
\usepackage{overpic}
\usepackage{subfigure}

\usepackage[algoruled,vlined]{algorithm2e}
\SetKwInOut{Input}{Input}
\SetKwInOut{Output}{Output}
\newcommand{\conv}{\operatorname{conv}}
\newcommand{\rank}{\operatorname{rank}}
\newcommand{\aff}{\operatorname{aff}}
\newcommand{\cP}{{\cal P}}
\newcommand{\cH}{{\cal H}}
\newcommand{\cF}{{\cal F}}
\newcommand{\cV}{{\cal V}}
\newcommand{\Pbar}{\bar{P}}%{\overline{P}}
\newcommand{\Fbar}{\bar{F}}
\newcommand{\R}{{\mathbb R}}
\newcommand{\Z}{{\mathbb Z}}
\newcommand\FACETENUMERATION{\textsc{FacetEnumeration}}
\newcommand\POLYTOPEVERIFICATION{\textsc{PolytopeVerification}}
\newcommand\POLYTOPECOMPLETENESS{\textsc{Completeness}}
\newcommand\POLYTOPEINCOMPLETENESS{\textsc{Incompleteness}}
\newcommand\POLYTOPECOMPLETENESSGEOMETRIC{\textsc{CompletenessG}}
\newcommand\POLYTOPECOMPLETENESSCOMBINATORIAL{\textsc{CompletenessC}}
\newcommand\FINDPULLINGFACET{\textsc{FindPullingFacet}}
\newcommand\ISPULLINGFACET{\textsc{IsPullingFacet}}
\newcommand\COMPLETENESSVIAHOMOLOGY{\textsc{Completeness\-Via\-Homology}}
\newtheorem{thm}{Theorem}[section]
\newtheorem{lem}[thm]{Lemma}
\newtheorem{cor}[thm]{Corollary}
\theoremstyle{definition}
\newtheorem{dfn}[thm]{Definition}
\newenvironment{problem}{\smallskip\parindent=0cm}{\smallskip}

\usepackage{amsmath}

\makeatletter
\renewcommand\@fnsymbol[1]{\ensuremath{\ifcase#1\or *\else +\fi}}
\makeatother

\title{Convex Hulls, Oracles, and Homology\footnote{An extended
    abstract version of this paper, ``Polytope verification by
    homology verification,'' has appeared in the Proceedings of
    EuroCG, Berlin, March 26--28, 2001, pp.~142--145.}}
\author{Michael Joswig\footnotemark[2] \qquad  G\"unter M. Ziegler\thanks{%
Partially supported by Deutsche Forschungs-Gemeinschaft (DFG)}\\
\small Inst.\ Mathematics, MA 6-2\\[-1.6mm]
\small TU Berlin, D-10623 Berlin, Germany\\[-1.6mm]
\small\url{{joswig,ziegler}@math.tu-berlin.de}}

\begin{document}
\date{\tiny January 10, 2003}
\maketitle

\begin{abstract}\noindent%
  This paper presents a new algorithm for the convex hull problem,
  which is based on a reduction to a combinatorial decision problem
  \POLYTOPECOMPLETENESSCOMBINATORIAL, which in turn can be solved by a
  simplicial homology computation.  Like other convex hull algorithms,
  our algorithm is polynomial (in the size of input plus output) for
  simplicial or simple input.  We show that the ``no''-case of
  \POLYTOPECOMPLETENESSCOMBINATORIAL{} has a certificate that can be
  checked in polynomial time (if integrity of the input is
  guaranteed).
\end{abstract}

\section{Introduction}
Every convex polytope $P\subset \R^d$ can be described as the convex
hull of a finite set $\cP$ of points or as the (bounded) set of
solutions of a finite system $\cH$ of linear equations and
inequalities \cite[Lect.~1]{GMZ}.  In view of the fundamental role
that polytopes play in Euclidean geometry and hence for any type of
geometric computing, the conversion between the two types of
representations, known as the \emph{convex hull problem}, is of key
interest.  It splits into two separate tasks.

The first task is the \emph{facet enumeration problem}: Given a finite
set of points \mbox{$\cP\subset\R^d$}, determine the combinatorial
structure of its boundary.  For this one does not want to explicitly
enumerate all the faces (the intersections of $P$ with supporting
hyperplanes), but one wants sparser date, namely to compute a minimal
representation of the convex hull $\conv(\cP)$ in terms of equations
and (facet-defining) inequalities.  Here the equations should describe
the affine hull $\aff(P)$, while the additional inequalities
correspond to the \emph{facets} (faces of codimension~$1$) of~$P$. If
$P$ is full-dimensional in~$\R^d$, then the facet-defining
inequalities are unique up to scaling.

The second task is the \emph{vertex enumeration problem}: Given a
finite system $\cH$ of linear (equations and) inequalities, and
provided that the set of solutions $P=\bigcap\cH$ is bounded, compute
the minimal set of points $\cP$ whose convex hull is~$P$. This minimal
set is unique; it consists of the \emph{vertices} ($0$-dimensional
faces) of~$P$.

The two tasks are dual to each other, via cone polarity.  Thus if an
LP-type oracle (an algorithm which for a system of inequalities
computes a solution, or for a set of points computes a separating
hyperplane, cf.~\cite{GLS}) is available, every algorithm for the
facet enumeration problem can also be used for vertex enumeration, and
vice versa.

Despite the great interest in the convex hull problem, and despite the
fact that a number of different strategies and algorithms have been
explored, implemented and analyzed in detail (see \cite{polyfaq}, as
well as Avis \cite{Avis} \cite{Avis2}, Fukuda \cite{Fuku} and Gawrilow
\& Joswig \cite{PolymakeCite} \cite{SoCG01:polymake} for
implementations), the problem can be considered ``solved'' neither in
theory, nor in practice.  If the dimension $d$ is fixed, Chazelle's
celebrated algorithm~\cite{MR94h:52026} gives an asymptotically
worst-case optimal (polynomial time) theoretical solution. Its
optimality is based on McMullen's ``Upper Bound Theorem'' \cite{McM}
on the maximal number of facets for a $d$-polytope with $n$ vertices.
However, for any given convex hull problem, the output may be small,
but it may also be much larger than the input --- indeed, it may be of
exponential size, if the dimension is not fixed.  This is very
relevant, since high-dimensional computations occur in a variety of
important applications.  Thus one is asking for a convex hull
algorithm whose running time is bounded by a polynomial in the size of
``input plus output''?  Such an algorithm would be called
\emph{output-sensitive}.  The analysis by Avis, Bremner and
Seidel~\cite{877.68119} shows that, unfortunately, none of the known
types of convex hull algorithms is output-sensitive.  These can
roughly be categorized as follows: Incremental and triangulation
producing (e.g., Chazelle's method), incremental without
triangulations (e.g., Fourier-Motzkin elimination
\cite[Lect.~1]{GMZ}), non-incremental (e.g., reverse
search~\cite{752.68082}).  Note that, by a result of
Bremner~\cite{990.25180}, only non-incremental methods can possibly be
output-sensitive.

The purpose of this paper is to describe a new (non-incremental)
convex hull algorithm, based on a completely different principle.  To
this end, we first present a (folklore) polynomial reduction of
\FACETENUMERATION{} to the decision problem \POLYTOPEVERIFICATION.
Then we further reduce to the \POLYTOPECOMPLETENESS{} problem: Is a
given description of a $d$-polytope by \emph{some} of its vertices and
\emph{some} of its facets complete, that is, are we given \emph{all}
the vertices and \emph{all} the facets? Looking at the convex hull
problem via its reduction to \POLYTOPEVERIFICATION{} or
\POLYTOPECOMPLETENESS{} automatically reveals its inherent self-dual
structure. It is an interesting feature that the
\POLYTOPECOMPLETENESS{} problem can be posed both with geometric input
data and as an entirely combinatorial problem
\POLYTOPECOMPLETENESSCOMBINATORIAL{}, where only the incidences
between vertices and facets are given.
 
Let us just mention here one recent occurrence of the combinatorial
completeness problem: McCarthy et al.\ \cite{MR1899891} describe a
situation where one wants to know whether a given inequality
description for a polytope is complete.  Moreover, the vertex
coordinates in some of their problems are necessarily non-rational, so
any coordinate-free/combinatorial approach is welcome. Unfortunately,
the most interesting case left ``open'' by McCarthy et al.\ (the
convex hull of the matrices corresponding to the Coxeter group $H_4$)
is a polytope completeness problem in dimension $d=16$ with $14{,}400$
vertices: From this data our method generates gigantic boundary
matrices that are plainly too large to process.

Also we have been informed by Samuel Fiorini (email, January 2002)
that he has successfully used a certificate for the ``no''-case of
\POLYTOPECOMPLETENESSCOMBINATORIAL{} that is similar to the one that
we describe in Section~\ref{sec:certificate}.

Our main contribution is an algorithm to attack the combinatorial
\POLYTOPECOMPLETENESSCOMBINATORIAL{} problem via deciding whether a
certain simplicial homology group of a certain abstract simplicial
complex vanishes or not. Moreover, we present a polynomially checkable
certificate for non-completeness, provided that the input is valid.
For the geometric version the validity of the input can be checked
easily.  Unfortunately, the complexity status for the homology
computation problem is open.  The best currently available strategy to
decide non-triviality of a (rational) homology group in question seems
to be to compute boundary matrices and perform Gauss elimination.
Since the boundary matrices in our algorithm can be exponentially
large, we do not obtain an output-sensitive method.  However, like
other methods (e.g., Avis' and Fukuda's reverse
search~\cite{752.68082} or Seidel's gift-wrapping
algorithm~\cite{Seid}) our algorithm is output-sensitive in the case
of simplicial polytopes.

\section{\FACETENUMERATION{} via \POLYTOPEVERIFICATION}

We start with a more formal description of the facet enumeration
problem:

\begin{problem}
\FACETENUMERATION$(e,\cP)$:\\
\Input{integer $e\ge0$; finite set of points $\cP\subset\R^e$.}
\Output{minimal description of $\conv(\cP)$ in terms of equations 
(for the affine hull of $\cP$) and
  inequalities (one for each facet of $\conv(\cP)$)}
\end{problem}

\noindent
It is known, cf.\ Avis, Bremner \& Seidel~\cite{877.68119},
Fukuda~\cite[Node~21]{polyfaq}, and Kaibel \& Pfetsch~\cite[Problems
1--3]{AGSS:KaibelPfetsch}, that \FACETENUMERATION{} has a polynomial
reduction to the polytope verification problem:

\begin{problem}
{\POLYTOPEVERIFICATION$(e,\cP,\cH)$}:\\
\Input{integer $e\ge0$; 
  finite set of points $\cP\subset\R^e$; 
  finite set $\cH$ of closed halfspaces in~$\R^e$}
\Output{answer \textbf{yes/no} to the question whether
  $\conv(\cP)=\bigcap\cH$}
\end{problem}

Freund and Orlin could show that a related problem, to decide whether
$\bigcap\cH\subseteq\conv(\cP)$, is co-NP-complete \cite{FreundOrlin}.

\section{\POLYTOPEVERIFICATION{} via \POLYTOPECOMPLETENESSGEOMETRIC}

Assuming that an LP-type oracle is available,
the \POLYTOPEVERIFICATION{} problem is polynomially equivalent to the
following \emph{geometric polytope completeness problem}:

\begin{problem}
{\POLYTOPECOMPLETENESSGEOMETRIC$(d,\cV,\cF)$}:\\
\Input{integer $d\ge0$; 
  finite set of points $\cV\subset\R^d$;
  finite set $\cF$ of closed halfspaces in~$\R^d$, such that\\
%\mbox{}\hspace{-18.3mm}such that
$\bullet$ $P:=\conv(\cV)$ is contained in $Q:=\bigcap\cF$\\
$\bullet$ $\dim P=\dim Q=d$\\
$\bullet$ every $v\in\cV$ defines a vertex of~$Q$\\
$\bullet$ every $F\in\cF$ defines a facet of~$P$}
\Output{answer \textbf{yes/no} to the question whether $P=Q$}
\end{problem}

\noindent
As in the case of \POLYTOPEVERIFICATION, the roles of vertices and
facets are interchangeable for \POLYTOPECOMPLETENESSGEOMETRIC.

We sketch the reduction of \POLYTOPEVERIFICATION{} to
\POLYTOPECOMPLETENESSGEOMETRIC. Given any input $(e,\cP,\cH)$ for
\POLYTOPEVERIFICATION, set $P:=\conv(\cP)$ and $Q:=\bigcap\cH$.
Employ Gaussian elimination to determine $\dim P$.  Verify 
%By computing $\#\cP\cdot\#\cH$ scalar products decide 
whether all the inequalities in~$\cH$ are valid for~$P$; if this is
not the case, then $P\not\subseteq Q$, so we output \textbf{no};
otherwise $P\subseteq Q$ is established. Now extract the set $\cH'$ of
all halfspaces from $\cH$ for which $P$ lies in the bounding
hyperplane, that is, all those inequalities which are tight
on~$\aff{P}$.  An LP-type oracle is sufficient, but also needed
\cite{GLS}, to check whether $\bigcap \cH'=\aff P$; if this is not the
case, then we know that $\dim Q>\dim P$, so we can output \textbf{no}.
Otherwise we proceed by restricting the input to $\aff P$, that is, we
deal with the situation where $P$ is full-dimensional.

Now remove from $\cH$ all the halfspaces which
% whose bounding hyperplanes 
do not determine facets of~$P$;
% that is, such that the intersection
%with $\cP$ has dimension smaller than $d-1$; 
this may be done using Gaussian elimination. (In the case $P=Q$, this
removal does not change $Q$; in the case $P\subset Q$, it may enlarge
$Q$.)  Similarly, we now remove from $\cP$ all those points which do
not arise as intersections of some bounding hyperplanes of halfspaces
in $\cH$; again this may be done via Gaussian elimination.
%This will remove from $\cP$ all points that are not vertices of~$P$. 
(In the case of $P=Q$, this removal does not change $P$; in the case
$P\subset Q$, we may loose vertices of~$P$, thus making $P$ smaller.)

Now we have prepared our input for \POLYTOPECOMPLETENESSGEOMETRIC.
Indeed, the first two conditions on the input are satisfied,
the other two are easily checked: % (again using Gaussian elimination):
If one of them fails, then output the answer \textbf{no}. 
%Otherwise feed the input to \POLYTOPECOMPLETENESSGEOMETRIC.
\hfill$\square$\par%\medskip 

\section{\POLYTOPECOMPLETENESSGEOMETRIC{} via
         \POLYTOPECOMPLETENESSCOMBINATORIAL}

The \emph{incidence matrix} of a polytope $P$ with
vertex set~$\cV$ and facet set~$\cF$ is defined to be
the matrix
\[
I_P\ \ :=\ \ (i_{Fv})_{F\in\cF,v\in\cV}\ \ \in\ \ \{0,1\}^{\cF\times\cV},
\]
where $i_{Fv}=1$ if vertex $v$ lies on the facet $F$ (that is, if
$v\in F$), and $i_{Fv}=0$ means that $v\notin F$.  This matrix is
well-defined up to permutation of rows and of columns, which
corresponds to reordering $\cV$ and $\cF$.  A \emph{minor} of a matrix
will refer to any submatrix obtained by possibly removing rows and/or
columns. A minor~$J$ of the incidence matrix~$I_P$ is \emph{complete}
if $J=I_P$.  Thus we arrive at the \emph{combinatorial polytope
  completeness problem}:

\begin{problem}
\POLYTOPECOMPLETENESSCOMBINATORIAL$(d,J)$:\\
\Input{integer $d\ge0$; incidence matrix minor $J$ of a $d$-polytope}
\Output{answer {\bf yes/no} to the question whether $J$ is complete}
\end{problem}

\noindent
It is not obvious that this problem is well defined.  However, from
Theorem~\ref{thm:homology} below it follows that there are no two
$d$-polytopes $P$ and~$P'$ such that a $0/1$-matrix~$J$ is both a
complete incidence matrix for~$P$ and an incomplete minor of an
incidence matrix for~$P'$.  (See also the related discussion
in~\cite{VertexFacetIncidences}.) It is clear that
\POLYTOPECOMPLETENESSGEOMETRIC{} has a polynomial reduction to
\POLYTOPECOMPLETENESSCOMBINATORIAL.

It is essential to have the dimension among the input parameters of
\POLYTOPECOMPLETENESSCOMBINATORIAL. This is demonstrated by the
following example \cite[p.~71]{GMZ}:
\[J_{\textrm{KM}}\ \ =\ \ \left(
\begin{array}{cccccccc}
1&1&1&1&0&0&0&0\\
1&1&0&0&1&1&0&0\\
1&0&0&1&1&0&0&1\\
0&1&1&0&0&1&1&0\\
0&0&1&1&1&1&0&0\\
0&0&0&0&1&1&1&1
\end{array}\right)
\]
We can identify $\cV = \{1,2,\dots,8\}$ and 
$\cF = \{1234, 1278, 1458, 2367, 3456, 5678\}$
with the sets of vertices and facets, respectively, of a
$3$-dimensional cube (in a suitable ``Klee-Minty'' %/``Gray code''
vertex numbering; see Figure~\ref{fig:pulling} below).  Consequently,
\POLYTOPECOMPLETENESSCOMBINATORIAL$(3,J_{\textrm{KM}})$ =
\textbf{yes}.  But we can also identify $\cV$ with the vertices of a
cyclic $4$-polytope $C_4(8)$.  Then each element in~$\cF$ corresponds
to a facet of~$C_4(8)$, according to Gale's evenness criterion. Hence
\POLYTOPECOMPLETENESSCOMBINATORIAL$(4,J_{\textrm{KM}})$ = \textbf{no},
since $C_4(8)$ has $20$~facets.

A more generic class of examples for which the dimension information
is needed arises from the prism construction: Let $P$ be an arbitrary
$d$-polytope and $P'=P\times[0,1]$ the prism over~$P$.  The facets of
$P'$ are $P\times\{0\}$, $P\times\{1\}$, and the products of facets
of~$P$ with the interval~$[0,1]$.  Call the latter facets of~$P'$
\emph{vertical}, and let $J_P$ be an incidence matrix of~$P$.  We have
\POLYTOPECOMPLETENESSCOMBINATORIAL$(d,J_P)$ = \textbf{yes}.  On the
other hand $J_P$ is also a minor of an incidence matrix of~$P'$, which
corresponds to the vertical facets and, say, the vertices in the
bottom facet $P\times\{0\}$. Therefore,
\POLYTOPECOMPLETENESSCOMBINATORIAL$(d+1,J_P)$ = \textbf{no}.

\section{\POLYTOPECOMPLETENESSCOMBINATORIAL{} via simplicial homology}

We will point out that \POLYTOPECOMPLETENESSCOMBINATORIAL{} has a
topological core.  The reader is referred to
Bj\"orner~\cite{Bj:TopMethods} for a survey of topological
combinatorics tools, and to Munkres~\cite{Munkres84:0} for a
presentation of simplicial homology. In the following we will use
reduced simplicial homology with coefficients in $\Z_2$.  One could
use any other commutative coefficient ring with unit, but $\Z_2$ is
the natural choice in terms of efficiency and simplicity.  We choose
non-reduced homology to simplify notation for the trivial case~$d=1$.

Let $J\in\{0,1\}^{\cF\times\cV}$ be an incidence matrix minor of some
polytope~$P$ with vertex set~$\cV'\supseteq\cV$ and facet set
$\cF'\supseteq\cF$.  Thus the columns of $J$ are in bijection with a
(partial) vertex set $\cV$ of~$P$. Each row of $J$ is the
characteristic vector of a subset of rows, i.e., of a subset of~$\cV$.
Thus in the following we interpret $J$ as a combinatorial encoding of
a system $\cF$ of (not necessarily distinct) subsets of~$\cV$, and
with slight abuse of notation we write $\cF\subseteq2^{\cV}$.  The
\emph{crosscut complex} of~$J$ is the simplicial complex
\[
\Gamma(J)\ \ :=\ \ 
\big(\cV,\ \bigcup\big\{2^F:F\in\cF\big\}\big),
\]
the simplicial complex of all sets of vertices that
are contained in \emph{some} facet in~$\cF$.

\begin{thm}\label{thm:homology}
  The incidence matrix minor~$J\in\{0,1\}^{\cF\times\cV}$ of a
  $d$-polytope is complete if and only if
  $\widetilde{H}_{d-1}(\Gamma(J);\Z_2)\neq0$.
\end{thm}

\begin{proof}
The set
\[
\Pi(P,J)\ \ :=\ \ \bigcup_{F\in\cF}\conv\{v\in\cV:v\in F\}
 \ \ \subseteq\ \ \partial P
\]
is a compact subset of the boundary of~$P$: For every ``given'' facet
$F$ of~$P$, it contains the convex hull of all ``given'' vertices.
Thus $\Pi(P,J)$ is a polyhedral complex, called a \emph{partial
  polytope}, covered by its convex (and hence contractible) cells
$\conv\{v\in\cV:v\in F\}$.  According to the nerve theorem
\cite{Bj:TopMethods}, the crosscut complex $\Gamma(J)$ has the same
homotopy type as the set~$\Pi(P,J)$.  In particular, the homology of
the set $\Pi(P,J)$ and of the crosscut complex coincide.  For an
example of the crosscut complex of a partial polytope see
Figure~\ref{fig:crosscut}.

In the \textbf{yes} case, if the sets of vertices and facets
both are complete, $\Pi(P,J)$ is the complete boundary of~$P$,
homeomorphic to $S^{d-1}$, so we have
$\widetilde{H}_{d-1}(\Gamma(\cF);\Z_2)\cong\Z_2$.

In the \textbf{no} case, if the vertex or the facet list is
incomplete, then $\Pi(P,J)$ is a proper subset of 
$\partial{P}$, 
which is a subcomplex of a suitable triangulation of $\partial P$,
so it cannot have $(d-1)$-dimensional homology.
\end{proof}

The complexity status of the problem to compute the rank of an
arbitrary homology group, or even to decide whether a certain homology
group vanishes, seems to be open; see Kaibel \&
Pfetsch~\cite[Problem~33]{AGSS:KaibelPfetsch}.  Thus currently our
best option is based on explicitly computing simplicial homology via
boundary matrices, as in Algorithm~\ref{alg:completenessviahomology}.

\begin{algorithm}[h]
\dontprintsemicolon
\caption{\COMPLETENESSVIAHOMOLOGY$(d,J)$%
\label{alg:completenessviahomology}}
\Input{integer $d\ge0$; an incidence matrix minor $J$ of a $d$-polytope}
\Output{answer {\bf yes/no} to the question whether $J$ is complete}
\BlankLine
\text{Generate $\Z_2$-boundary matrices $\partial_d$ and
$\partial_{d-1}$ for $\Gamma(J)$}\;
\eIf{$\dim_{\Z_2}\ker\partial_{d-1}>\rank_{\Z_2}\partial_d$}{
  \Return \textbf{yes}}{\Return \textbf{no}}
\end{algorithm}

To estimate the costs of this computation, suppose that $n=|\cV|$,
$m=|\cF|$, and that the maximum cardinality of any facet equals~$s$.
Thus $J\in\{0,1\}^{m\times n}$, and every row of $J$ contains at 
most $s$ ones.
Then the size of the relevant boundary matrices is bounded from above
by $\binom{s}{d+1}m\times\binom{s}{d}m$ and
$\binom{s}{d}m\times\binom{s}{d-1}m$, respectively. We use Gaussian
elimination over~$\Z_2$ to compute the rank and the corank,
respectively.

\begin{cor}
  The algorithm \COMPLETENESSVIAHOMOLOGY$(d,J)$ has a polynomial
  running time if $s$ is bounded by $d+c$, for an absolute constant
  $c\ge0$.
\end{cor}

The latter case is, in fact, interesting: A $d$-polytope is 
\emph{simplicial} if each proper face is a simplex or, equivalently,
each facet contains exactly $d$ vertices.  We infer that the running
time of \COMPLETENESSVIAHOMOLOGY{} for simplicial polytopes is bounded
by~$O(dm^3)$.

It has been observed by Bremner, Fukuda \& Marzetta~\cite{910.68217}
that \FACETENUMERATION{} for a polytope~$P$ is polynomially equivalent
to \FACETENUMERATION{} for the dual polytope~$P^*$.  Using our
techniques, a similar result can be obtained directly.  If $I$ is an
incidence matrix for~$P$, then the transposed matrix~$I^{\text{tr}}$
is an incidence matrix for~$P^*$.  Any minor $J$ of~$I$ is complete if
and only if its transpose is a complete minor of~$I^{\text{tr}}$. This
leads to the following modification of our algorithm.  While $s$ was
defined above as the maximal row size of the input incidence matrix
minor, define 
\[
s'\ :=\ \min\{\text{maximal row size},\text{maximal column size}\}.
\]
Thus we modify our algorithm: It should first compare the sizes of the
primal and the dual problem, and then perform the (reduced) homology
computation for the smaller problem.
%\begin{cor}
  The modified algorithm \COMPLETENESSVIAHOMOLOGY$(d,J)$ has
  polynomial running time if $s'$ is bounded by ``$d$ plus a constant.''
%\end{cor}
In particular, this yields an $O(d(n+m)^3)$-algorithm for the
\POLYTOPECOMPLETENESSCOMBINATORIAL{} problem specialized to polytopes
which are simplicial or \emph{simple}, that is, dual to a simplicial
polytope.

We note, however, that these running times are neither optimal nor the
best available: The reverse search algorithm of Avis and
Fukuda~\cite{752.68082} computes the convex hull (and thereby solves
\POLYTOPECOMPLETENESSGEOMETRIC) of a simplicial polytope in~$O(dnm)$
steps.

\section{A Certificate for \POLYTOPEINCOMPLETENESS}
\label{sec:certificate}

\begin{figure}[hb]
  \begin{center}
    \subfigure[The ($3$-dimensional) crosscut complex of some partial
    $3$-cube~$C$.  The two quadrangle faces of~$C$ yield tetrahedra in
    $\Gamma(C)$, which are displayed almost
    flat.\label{fig:crosscut}]{\includegraphics[width=7.5cm]{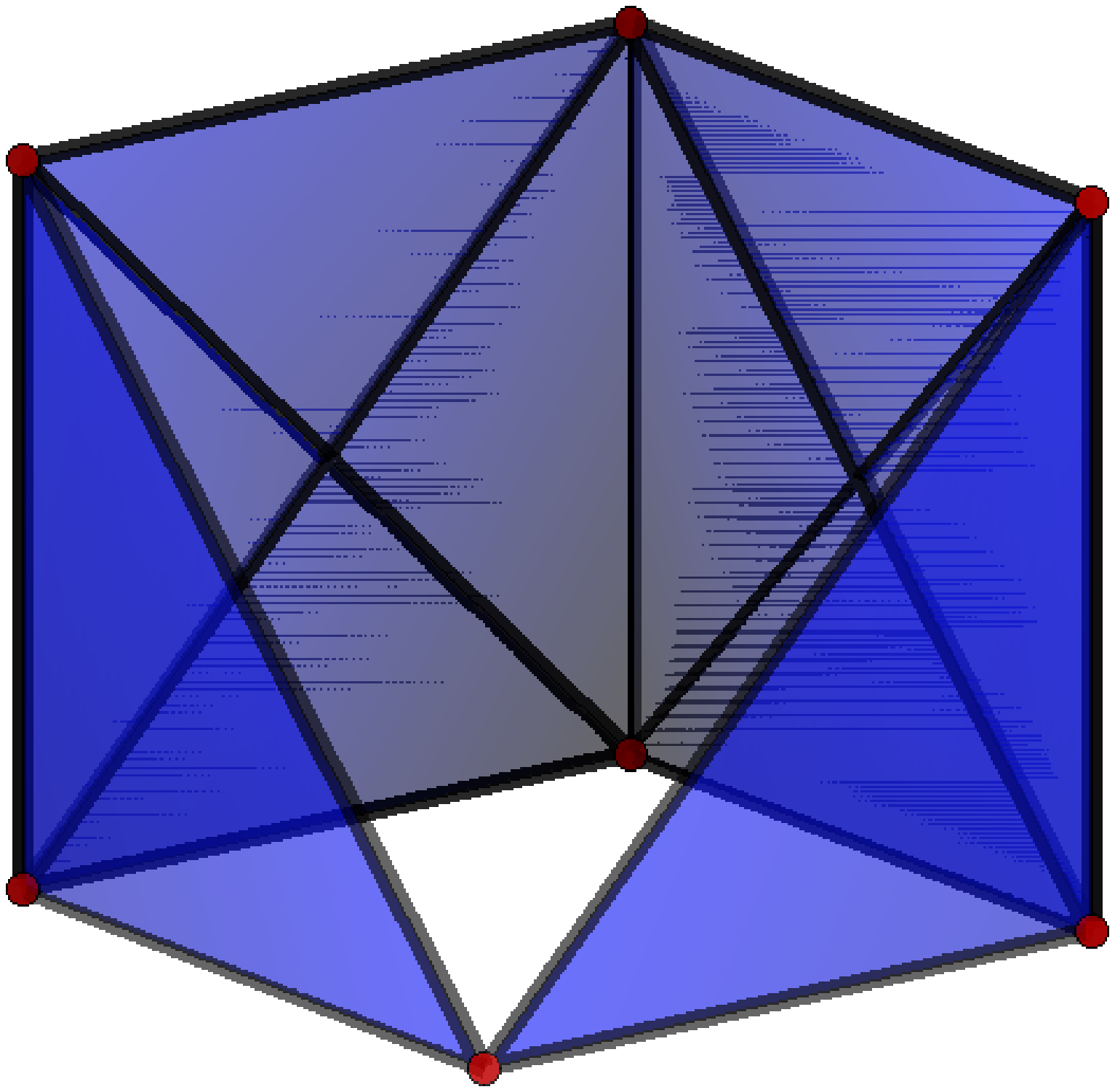}}
    \subfigure[The pulling triangulation of the boundary of a $3$-cube
    with respect to a ``Klee-Minty'' vertex ordering.  The facet
    $\{1,7,8\}$ of the triangulation corresponds to the flag
      $\{8\}\subset\{7,8\}\subset\{1,2,7,8\}$ of the cube.\label{fig:pulling}]{%
      \begin{overpic}[width=7.5cm]{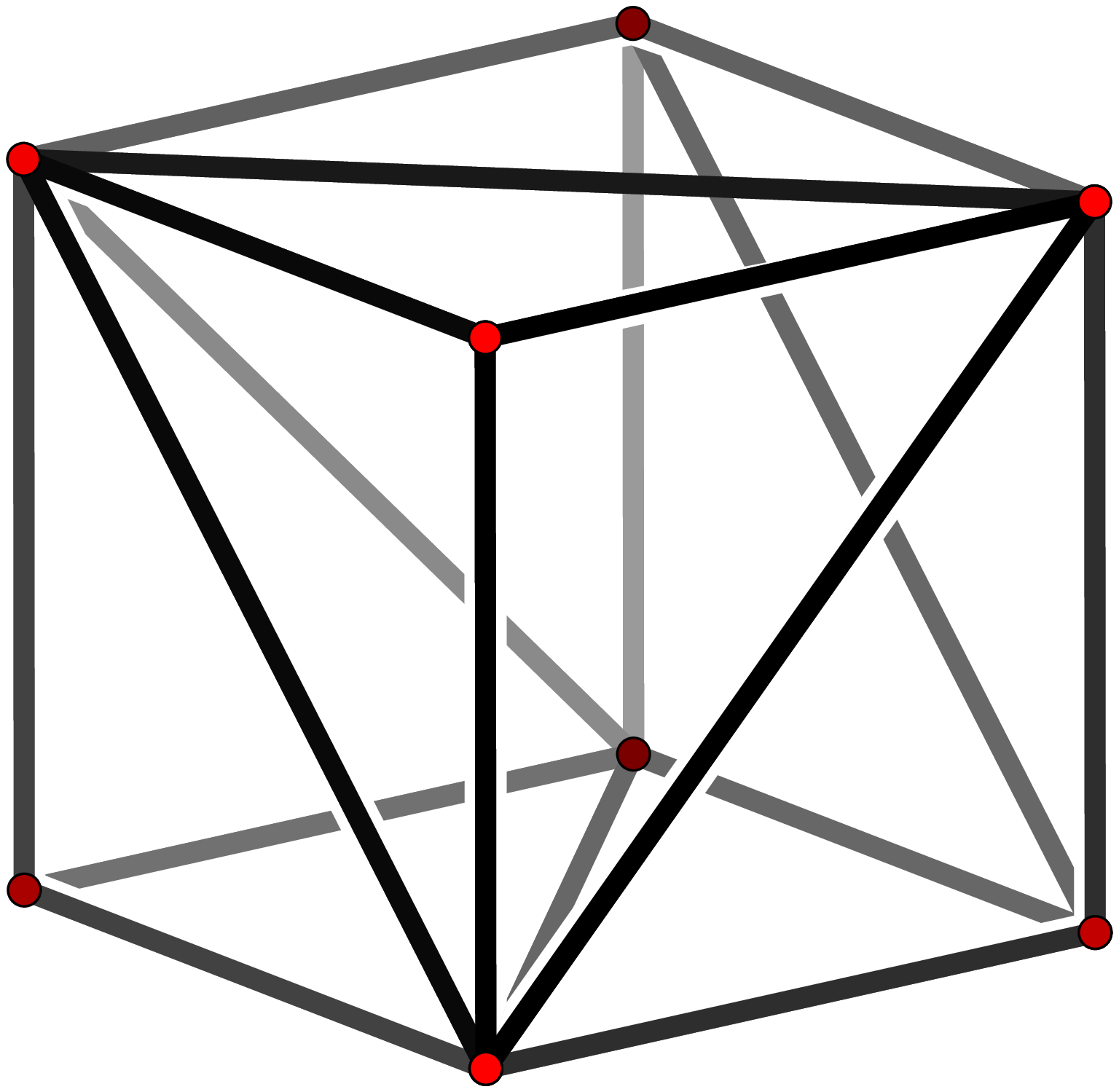}
        \put(44,-1){$1$}
        \put(88,10){$2$}
        \put(57,30){$3$}
        \put( 9,13){$4$}
        \put( 9,68){$5$}
        \put(57,77){$6$}
        \put(88,65){$7$}
        \put(44,58){$8$}
      \end{overpic}}

    \caption{Crosscut complex and pulling triangulation.}
  \end{center}
\end{figure}

Let $P$ be a $d$-polytope with ordered vertex set
$\cV'=\{v_1,\dots,v_n\}$ and facet set~$\cF'$.  Inductively, define a
sequence $\Delta_0,\dots,\Delta_n$ of polytopal subdivisions of the
boundary complex~$\partial P$: Set $\Delta_0:=\partial P$.  In order
to obtain $\Delta_k$ replace each facet~$F$ of~$\Delta_{k-1}$ which
contains~$v_k$ by the set of cones with apex~$v_k$ over those facets
of~$F$ which do not contain~$v_k$.  The final subdivision is a
triangulation $\Delta(P):=\Delta_m$ of~$\partial P$, the \emph{pulling
  triangulation} \cite{HDCG:Lee} with respect to the chosen ordering
of~$\cV$.  For an example of a pulling triangulation see
Figure~\ref{fig:pulling}.

The pulling triangulation of $\partial P$ has several nice properties
(not shared, for example, by the ``placing triangulation'') which may
be exploited for our purposes.  First, its combinatorics is determined
by the combinatorics of $P$; see below.  Furthermore, if we use a
linear ordering of the vertex set $\cV'$ in which the vertices
in~$\cV$ come first, then the corresponding pulling triangulation of
the boundary of~$P$ contains $\Pi(P,J)$ as a subcomplex.

Let us now identify the vertex set~$\cV'$ with the set
$[n]=\{1,\dots,n\}$ and each facet~$F\in\cF'$ with the subset of~$[n]$
that corresponds to the vertices contained in~$F$.  Thus any
triangulation of~$\partial P$ is encoded by a collection of
$d$-subsets of~$[n]$, that is, to a subset of $\binom{[n]}{d}$.  We
write $\{v_1,\dots,v_d\}_<$ for a $d$-subset of~$[n]$ with
$v_1<v_2<\dots<v_d$.

\begin{lem}\label{lem:pulling}
Let $P$ be a $d$-polytope whose vertex set is labeled by~$[n]$.\\
Then a set $\{v_1,\dots,v_d\}_<\in\binom{[n]}{d}$ 
corresponds to a facet of the pulling
triangulation of~$\partial P$ (with respect to the chosen vertex labeling)
if and only if there is a complete flag of faces
\[
\emptyset\ \subset\ 
      G_0\ \subset\ 
      G_1\ \subset\ \ \dots\ \ \subset\ 
  G_{d-1}\ \subset\ P,
\]
such that $v_i$ is the smallest vertex in~$G_{d-i}$ for $1\le i\le d$,
that is, if there are facets $F_1,\dots,F_d$ of~$P$ such that 
\[ 
v_i\ =\ \min(F_1\cap\ldots\cap F_i)
\]
for~$1\le i\le d$.  
\end{lem}

\begin{proof}
Every pulling facet $\{v_1,\dots,v_d\}_<$ lies in a 
facet $F_1=G_{d-1}$ of~$P$, with $v_1=\min G_{d-1}$.
It is a cone with apex $v_1$ and base $G_{d-2}\subset G_{d-1}$.
The existence of the rest of the maximal flag $(G_i)_{0\le i<d}$
follows recursively. Given the flag, the existence of the 
facets $F_1,\dots,F_d$ follows \cite[Lect.~2]{GMZ}.
%Not every maximal flag induces a pulling facet. 
Given a complete flag, the corresponding
sequence of facets $F_i$ is uniquely determined if $P$ is
simple, but not in general.
\end{proof}

If we have an arbitrary incidence matrix minor~$J$ of a
$d$-polytope~$P$, then we can read the combinatorial characterization
of the pulling triangulation from Lemma~\ref{lem:pulling} as the
definition of a complex that coincides with the pulling triangulation
of~$\partial{P}$  in case $J$ is complete, but is well-defined
in general:

\begin{dfn}
Given an integer $d>0$ and a 
$0/1$-matrix $J\in\{0,1\}^{m\times n}$, which we 
interpret as the incidence matrix of a set system
$\cF\subseteq2^{[n]}$, the \emph{pulling complex} of $d$~and~$J$~is 
\begin{eqnarray*}
\Delta(d,J)\ \ :=\ \ 
\Big\{\{v_1,\dots,v_d\}_< \in\tbinom{[n]}d & :& 
\textrm{there are } \Fbar_1,\dots,\Fbar_d\in\cF\textrm{ such that } 
\\[-3.6mm]
&&v_i=\min(\Fbar_1\cap\ldots\cap \Fbar_i)
\textrm{ for }1\le i\le d
\Big\}.
\end{eqnarray*}
\end{dfn}

%It seems important to ``interpret'' $\Delta(d,J)$
%also in the case where $J$ is an incomplete minor.

\begin{lem}
  Let $P$ be a $d$-dimensional polytope with vertex set $\cV'$ and
  facet set $\cF'$, and let $J$ be a incidence matrix minor
  corresponding to subsets $\cV\subseteq\cV'$ and $\cF\subseteq\cF'$.
  Let $\Pbar\subseteq P$ be the convex hull of the vertices in~$\cV$.
  Fix a linear ordering on the vertex set $\cV'$ such that the
  vertices in~$\cV$ come first.
  
  Then the simplicial complex $\Delta(d,J)$ is a subcomplex of
  $\Delta(P)$ as well as of $\Delta(\Pbar)$.  In particular,
  $\Delta(d,J)$ is a proper subcomplex of $\Delta(P)$, unless the
  minor $J$ is complete, $J=I_P$.  In the incomplete case
  $\Delta(d,J)$ may even be empty.
\end{lem}

\begin{proof}
  Let $\{v_1,\dots,v_d\}_< \in\Delta(d,J)$, then there are
  $\Fbar_1,\dots,\Fbar_d\in\cF$ such that
  $v_i=\min(\Fbar_1\cap\ldots\cap\Fbar_i)$.  Now since $J$ is an
  incidence matrix minor of~$P$, there are facets
  $F_i\supseteq\Fbar_i$ of~$P$, and by the assumption on the vertex
  ordering the vertices in $\Fbar_i$ come first, so
  $\min(\Fbar_1\cap\ldots\cap\Fbar_i)=\min(F_1\cap\ldots\cap F_i)$,
  which yields $\{v_1,\dots,v_d\}_< \in\Delta(P)$.
  
  Now $\Pbar=\conv(\cV)$, and the $\Fbar_i=F_i\cap\cV$ are vertex sets
  of faces (not necessarily facets) of~$\Pbar$.  If the vertices
  $v_i=\min(\Fbar_1\cap\ldots\cap\Fbar_i)$ are distinct, then the
  faces $\Fbar_1\cap\ldots\cap\Fbar_i$ form a complete flag in the
  face lattice of~$\Pbar$, and thus $\{v_1,\dots,v_d\}_<
  \in\Delta(\Pbar)$, by Lemma~\ref{lem:pulling}.
\end{proof}

In particular, $\Delta(d,J)$ triangulates a subset of the complex
$\Pi(P,J)$ that appears in the proof of Theorem~\ref{thm:homology}.

Now we present a polynomially-checkable certificate for the case that
$J$ is incomplete. Note, however, that this result does \emph{not}
prove that \POLYTOPECOMPLETENESSCOMBINATORIAL{} is in co-NP: We are
not able to check (in polynomial time) whether the input is valid,
that is, whether $J$ is actually an incidence matrix minor of some
$d$-polytope.

\begin{thm}\label{thm:certify}%
  Any {\bf no} instance of the problem
  \POLYTOPECOMPLETENESSCOMBINATORIAL$(d,J)$ has a certificate that can
  be verified in polynomial time.
\end{thm}

\begin{proof}
  The minor $J$ is incomplete if and only if the pulling complex
  $\Delta(d,J)$ is not a complete triangulation of a $d$-polytope
  boundary.  Two cases arise. The first one is if
  $\Delta(d,J)=\emptyset$, in which case Algorithm~\ref{alg:find}
  described below will certify in polynomial time that $J$ is not
  complete.
  
  The second case is if $\Delta(d,J)$ is non-empty but incomplete.  In
  this case (since the dual graph of the pulling triangulation
  $\Delta(P)$ is connected) there is a facet $\{v_1,\dots,v_d\}\in
  \Delta(d,J)$ together with an index $i$ such that there is no second
  facet of $\Delta(d,J)$ that contains $\{v_1,\dots,v_d\}\setminus
  \{v_i\}$.  In this situation our certificate is the set
  $\{v_1,\dots,v_d\}\setminus\{v_i\}$.  Calling \ISPULLINGFACET{} for
  every $d$-subset of~$[n]$ which contains the certificate, this
  certificate can be verified in polynomial time, since there are
  $n-d+1$ of these subsets.
\end{proof}

Now we proceed by describing the two subroutines needed for
Theorem~\ref{thm:certify}. The first one is Algorithm~\ref{alg:find}:
Given an incidence matrix minor~$J$ it either finds a facet
of~$\Delta(d,J)$ in polynomial time or it detects that $J$ is
incomplete.  The correctness follows from Lemma~\ref{lem:pulling}.
Our specific formulation of the algorithm
%Since at the beginning of the first round of
%the for-loop $S=[n]$ the facet $F_1$, constructed in line~\ref{line}
%is excluded to contain the vertex $1=\min[n]$.  Therefore the
%algorithm never 
produces a pulling triangulation facet which does not
contain~$1$:  This restriction does not hurt, since $\Delta(d,J)$
must contain such a facet if $J$ is complete.

\begin{algorithm}[h]
\dontprintsemicolon
\caption{\FINDPULLINGFACET$(d,J)$\label{alg:find}}
\Input{incidence matrix minor $J\in\{0,1\}^{m\times n}$ of a
  $d$-polytope; $d$-tuple\\ 
  \mbox{$\{v_1,\dots,v_d\}_<\in\binom{[n]}{d}$}}
%  $(d,J)$ as above}
\Output{a facet $\{v_1,\dots,v_d\}\in\Delta(d,J)$, 
or \textbf{incomplete}}
\BlankLine
$S\leftarrow[n]$\;
\For{$i\leftarrow 1$ to $d$}{
%\nl\label{line}% 
  $F_i\leftarrow$ any $F\in\cF$ such that $\min S\notin F$, 
  $F\cap S\ne\emptyset$, and $|F\cap S|$ is maximal\;
  \If{no such facet exists}{
    \Return \textbf{incomplete}
    }
  $S\leftarrow S\cap F_i$\;
  $v_i\leftarrow\min S$\;
  }
\Return $\{v_1,\dots,v_d\}_<$
\end{algorithm}

Our second subroutine, Algorithm~\ref{alg:ispullingfacet}, checks
whether a given set of $d$~vertices is a facet of the pulling
complex~$\Delta(d,J)$ or not.  Its correctness again follows from the
characterization in Lemma~\ref{lem:pulling}. Its running time is
bounded by $O(d(n+m))$.

\begin{algorithm}[ht]
\dontprintsemicolon
\caption{\ISPULLINGFACET$(d,J,\{v_1,\dots,v_d\}_<)$%
\label{alg:ispullingfacet}}
\Input{$(d,J)$ as above}
%  incidence matrix minor $J\in\{0,1\}^{m\times n}$ of some
%  $d$-polytope; $\{v_1,\dots,v_d\}_<\in\binom{[n]}{d}$}
\Output{answer \textbf{yes/no} to the question whether
  $\{v_1,\dots,v_d\}\in\Delta(d,J)$}
\BlankLine
\For{$i\leftarrow d$ downto $1$}{
  compute the set $\cF_i$ of all facets (i.e., rows of $J$) 
  that contain $\{v_i,\dots,v_d\}$\;}
\For{$i\leftarrow 1$ to $d$}{
  $F_i\leftarrow$ any $F\in\cF_i$ with 
  $v_i=\min(F_1\cap\ldots\cap F_{i-1}\cap F)$\;
  \If{no such $F$ exists}{
    \Return \textbf{no}\;
    }
  }
\Return \textbf{yes}
\end{algorithm}

%If $J$ is complete, that is, $\Delta(d,J)$ is a pulling triangulation
%of~$\partial P$, then the dual graph of the triangulation is
%connected.  Thus by starting at an arbitrary facet of~$\Delta(d,J)$,
%obtained from \FINDPULLINGFACET$(d,J)$, we can compute all facets
%of~$\Delta(d,J)$ as follows: Proceed, e.g., by breadth first search
%through the dual graph; each facet has precisely $d$ neighbors; each
%neighbor of a facet~$\sigma\in\Delta(d,J)$ differs from any neighbor
%in one vertex; all candidates can be considered and verified by
%calling \ISPULLINGFACET.  This yields algorithms for computing a
%pulling triangulation of the boundary of a polytope which are
%polynomial in the size of the output (that is, the triangulation).
%However, there are polytopes (e.g., products of simplices) which do
%not have triangulations of size polynomial in~$m$ and~$n$; see
%Lee~\cite{HDCG:Lee}.

%Let us also note at this point that the construction of a pulling
%triangulation also leads to algorithms for the \FACETENUMERATION{}
%problem.  These belong to the class of ``triangulation producing
%algorithms'' in the sense of Avis, Bremner \& Seidel~\cite{877.68119}.
%In particular, such algorithms cannot be polynomial in the size of
%``input plus output.''

We close our discussion with a pointer to a specific special case: It
would be interesting to know whether \POLYTOPECOMPLETENESS$(d,J)$
has a polynomial time solution for the very special case where $J$
has all columns and lacks at most one row.

\subsection*{Acknowledgements}

We are grateful to Volker Kaibel, Marc E. Pfetsch and Mark de
Longueville for helpful comments.  Moreover, the first author is
indebted to G\"unter Rote for an enlightening discussion on the
subject.

\begin{small}
\providecommand{\bysame}{\leavevmode\hbox to3em{\hrulefill}\thinspace}

\end{small}
\bibliographystyle{siam}
%\bibliography{../../bib/POLYref}%mic,top,poly,geo,math,soft,AGSS}
\end{document}